\documentclass[reqno]{amsart}
\usepackage[dvips]{graphicx}

\newtheorem{thm}{Theorem}

\newtheorem{dfn}{Definition}

\newtheorem{prp}{Proposition}

\newtheorem{lma}{Lemma}

\newtheorem{qst}{Question}

\newtheorem{rem}{Remark}
\newenvironment{rmk}{\begin{rem}\rm}{\end{rem}}

\newenvironment{pf}{\begin{proof}}{\end{proof}}

\newcommand{\tR}{{\mathbb{R}}}

\newcommand{\tZ}{{\mathbb{Z}}}
\newcommand{\tS}{{\mathbb{S}}}

\newcommand{\la}{\langle}
\newcommand{\ra}{\rangle}

\newcommand{\bi}{\operatorname{bridge}}
\newcommand{\ba}{\operatorname{braid}}
\newcommand{\hc}{\kappa_{{\rm hol}}}
\newcommand{\hct}{(\kappa+\tau)_{{\rm hol}}}
\newcommand{\E}{{\mathcal E}}
\newcommand{\Eh}{{\mathcal E}_{{\rm hol}}}

\begin{document}

\title
{Total curvatures of holonomic links}
\author{Tobias Ekholm \and Ola Weistrand}
\address{Dept. of math. Uppsala University,S-751 06 Uppsala, Sweden}
\email{tobias@math.uu.se}  

\begin{abstract}
A differential geometric characterization of the braid-index of a link 
is found. After multiplication by $2\pi$, it equals the infimum
of the sum of total curvature and total absolute torsion over
holonomic representatives of the link.

Upper and lower bounds for the infimum of total curvature over
holonomic representatives of a link are given in terms of its braid-
and bridge-index. Examples showing that these bounds are sharp are
constructed.
\end{abstract}

\maketitle

\section{Introduction}
Let $C_k$ denote the disjoint union of $k$ circles. An isotopy class
of embeddings $C_k\to\tR^3$ will be called a {\em ($k$-component) link}.
A $1$-component link will be called a {\em knot}.

A collection of loops $c\colon C_k\to\tR^3$ is {\em holonomic} if it
arises as the $2$-jet extension of a function $f\colon C_k\to\tR$. That 
is, in coordinates $(x_1,x_2,x_3)$ on $\tR^3$,
$c(t)=\left(f(t),f'(t),f''(t)\right)$, $t\in C_k$. Vassiliev ~\cite{V}
introduced holonomic loops in knot theory and proved that any tame
link has a holonomic representative. Later, Birman and Wrinkle
~\cite{BW} proved that if two holonomic embeddings $C_k\to\tR^3$ are
isotopic then they are isotopic through holonomic embeddings.  

The {\em curvature} $\kappa(L)$ and {\em curvature-torsion}
$(\kappa+\tau)(L)$ of a link $L$ are the infima of total curvature
and of the sum of total curvature and total absolute torsion, 
respectively, over all embeddings representing $L$. 
These invariants were defined by Milnor ~\cite{M1}, ~\cite{M2}. He
proved that if $L$ is any link then $\kappa(L)=2\pi\bi(L)$, where  
$\bi(L)$, the {\em bridge-index} of $L$, equals the minimal number of
over-arcs in a link-diagram representing $L$, and  
$(\kappa+\tau)(L)$ is an integral multiple of $2\pi$. 

Further results on the curvature-torsion invariant were found by 
Honma and Saeki ~\cite{HS}. They proved if $L$ is any link then 
$(\kappa + \tau)(L)\le2\kappa(L)-2\pi$ and showed that 
the difference between $2\pi\ba(L)$, where $\ba(L)$, the {\em
braid-index} of $L$, equals the minimal number of strands in a closed braid
representing $L$, and $(\kappa+\tau)(L)$ may be strictly positive (and
actually arbitrarily large).

In this paper, we shall study the {\em holonomic curvature} $\hc$ and
{\em holonomic curva-ture-torsion} $\hct$. 
For a link $L$, $\hc(L)$ and
$\hct(L)$ are the infima of total curvature and of
the sum of total curvature and total absolute torsion,
respectively, over all {\em holonomic} embeddings representing
$L$ (see Section ~\ref{secctoh}). 
\begin{thm}\label{thmhct}
For any link $L$,
$$
\hct(L)=2\pi\ba(L).
$$
\end{thm} 
\noindent
Theorem ~\ref{thmhct} is proved in Section ~\ref{pfthmhct}. It gives a
differential geometric characterization of the braid-index.
\begin{thm}\label{thmhc}
For any link $L$,
\begin{equation}\label{ineqhc}
A(\ba(L)-\bi(L))+2\pi\bi(L)\le\hc(L)\le2\pi\ba(L),
\end{equation}
where $2A$ is the area on the unit sphere
$\tS^2=\{(x,y,z)\in\tR^3\colon x^2+y^2+z^2=1\}$ of the region defined
by the inequality $y^2-4xz>0$.
\end{thm}
\noindent
Theorem ~\ref{thmhc} is proved in Section ~\ref{pfthmhc}. The constant
$A=(1.29\dots)\pi$ can be expressed in terms of elliptic
integrals.    

For links $L$ with $\bi(L)=\ba(L)$ the inequalities ~\eqref{ineqhc}
are equalities. Hence, the upper bound in ~\eqref{ineqhc} is best
possible. So is the lower bound:
\begin{prp}\label{prpex}
For any  $m\ge 2$, $j\ge 1$ there exists a knot $K(m,j)$ such that
$\bi(K(m,j))=m$, $\ba(K(m,j))=m+j$, and $\hc(K)=Aj+2\pi m$.
\end{prp}
\noindent
Proposition ~\ref{prpex} is proved in Section ~\ref{pfprpex}. Knots
with properties as stated in Proposition ~\ref{prpex} are defined in
Section ~\ref{K(m,j)}. 

Theorem ~\ref{thmhc} gives rise to the following question:
\begin{qst}\label{qst}
Is it true that 
$$
A(\ba(L)-\bi(L))+2\pi\bi(L)=\hc(L),
$$  
for all links $L$?
\end{qst}
\noindent
An affirmative answer would give a differential
geometric characterization of the difference between braid- and
bridge-index. A negative answer would give a new link invariant:
$\hc-2\pi\bi$. 

The space $\Eh$ of holonomic embeddings $C_k\to\tR^3$ is a subspace of
the space $\E$ of all embeddings . The result of Birman and
Wrinkel mentioned above implies that the inclusion $\Eh\subset\E$
induces an isomorphism on $\pi_0$. Theorems ~\ref{thmhct} and
~\ref{thmhc} shows that there are (in a sense large) open subsets of
$\E$ which do not intersect $\Eh$. 

\section{Holonomic curvature and curvature-torsion}\label{secctoh}
\subsection{Holonomic curvature}
If $c$ is a continuous closed space curve we denote its {\em total
curvature}, as defined in ~\cite{M1}, p.251, by 
$\int_c\kappa\,ds$.
\begin{lma}\label{avtc}
Let $c$ be a continuous closed space curve. Then
\begin{equation}\label{deftc}
\int_c\kappa\,ds=\frac12\int_{\tS^2}\mu(c,v)\,dA(v),
\end{equation}
where $dA$ is the area form on the unit sphere $\tS^2\subset\tR^3$ and
$\mu(c,v)$ is the number of local maxima of the function 
$t\mapsto \la c(t),v\ra$. (The symbol $\la\, ,\ra$ denotes the
standard inner-product on $\tR^3$.)  
\end{lma}
\begin{pf}
This is Theorem 3.1 in ~\cite{M1}.
\end{pf}
\begin{rmk}\label{rmkext}
Equation ~\eqref{deftc} could be taken as definition of total
curvature. 
\end{rmk}
\begin{rmk}
If $c$ is a space curve with continuous unit tangent vector $e_1$ then 
the total curvature of $c$ equals the length on the unit sphere
$\tS^2\subset\tR^3$ of the curve traced out by $e_1$. 

If $c$ is a twice continuously differentiable curve with
everywhere defined curvature function $\kappa$ then the total
curvature of $c$ equals the integral over $c$ of
$\kappa\,ds$, where $ds$ is the arclength-element. 
\end{rmk}
\begin{dfn}
The {\em holonomic curvature} of a link $L$ is the infimum of
$\int_c\kappa\,ds$ over all holonomic representatives $c$ of $L$. 
\end{dfn}
\subsection{Holonomic curvature-torsion}
A three times continuously differentiable curve $c\colon[a,b]\to\tR^3$ 
is called {\em non-degenerate} if $c'(t)$ and 
$c''(t)$ are linearly independent for every $t\in[a,b]$. If $c$ is
non-degenerate then its curvature and torsion functions $\kappa$ and $\tau$ 
are defined and $\kappa>0$.
\begin{dfn}
The {\em sum of total curvature and total absolute torsion} of a
non-degenerate curve $c$ is  
$$
\int_c(\kappa+|\tau|)\,ds,
$$
where $ds$ is the arclength-element.
\end{dfn}

We shall be concerned with holonomic curves and so would like to
describe the set of functions $f\colon C_k\to\tR$ with $2$-jet
extensions which are non-degenerate curves. 
\begin{dfn}\label{torgen}
A function $f\colon C_k\to\tR$ is {\em torsion-generic} if its
associated holonomic curve is non-degenerate.  
\end{dfn}

\begin{prp}\label{cod2}
The set of non-torsion-generic functions $C_k\to\tR^3$ has
codimension $2$ in the space of all smooth functions
$C_k\to\tR^3$. (That is, non-torsion-generic functions can be avoided
in generic $1$-parameter families of functions and in generic
$2$-parameter families they appear at isolated points.)
\end{prp}
\begin{pf}
Let $(t,x_0,x_1,x_2,x_3,x_4)\in S^1\times\tR^4$ be coordinates on the
jet-space $J^4(S^1,\tR)$ and let $N$ be the locus of the equations 
$$
x_1x_3-x_2^2=0,\quad x_1x_4-x_2x_3=0,\quad x_2x_4-x_3^2=0.
$$
Then $N$ is an algebraic subvariety of $J^4(S^1\times\tR^4)$. The
non-singular part of $N$ is a submanifold of codimension $2$ and the
singular part of $N$ equals $V_1\cup V_2$ where $V_1$ is the locus of
$x_2=x_3=x_4=0$ and $V_2$ the locus of $x_1=x_2=x_3=0$.

A function $f\colon S^1\to\tR$ is non-torsion-generic if and only if 
its $4$-jet extension $j^4f$ satisfies $j^4f(S^1)\cap N\ne\emptyset$.
Applying the jet-transversality theorem, the proposition follows.    
\end{pf}
\begin{dfn}
The {\em holonomic curvature-torsion} of a link $L$ is the infimum of
$\int_c(\kappa+|\tau|)\,ds$ over all torsion-generic holonomic
representatives $c$ of $L$.  
\end{dfn}
\section{Proof of Theorem ~\ref{thmhct}}\label{pfthmhct}
Let $(x_1,x_2,x_3)$ be coordinates on $\tR^3$ and let $j\ge 1$ be an
integer. Let $f\colon C_j\to\tR$ 
be a {\em generic} (as defined in ~\cite{V}, Proposition 1, (i), (a))
and torsion-generic function such that the embedding 
$c\colon C_j\to\tR^3$, 
$c(t)=\left(f(t),f'(t),f''(t)\right)$, $t\in C_j$   
represents the link $L$. 

Let $n$ be the number of local maxima of $f$. By performing a sequence
of holonomic second Reidemeister moves (more precisely, moves as shown 
in ~\cite{V}, Fig. 1e) we can assure that any
local maximum of $f$ is larger than the largest local minimum. The
projection of the holonomic curve of $f$ onto the $x_1x_2$-plane, endowed
with over/under-information, is then a closed braid on $n$ strands. Thus,
$n\ge\ba(L)$. Since the linking number of $c$ and the $x_1$-axis is $-n$,
it follows from ~\cite{M2} Theorem 3 that 

$$
\int_c(\kappa + |\tau|)\,ds\ge 2\pi n.
$$
Noting that the set of generic functions is open and dense in the
space of all functions $C_k\to\tR$, we conclude that $\hct(L)\ge2\pi\ba(L)$.

Theorem 1 in ~\cite{BW} implies that $L$ has a holonomic representative
such that its diagram $D$ in the
$x_1x_2$-plane is the diagram of a closed braid on $\ba(L)$
strands. (The diagram $D$ has all its negative crossings in the upper
half-plane ($x_2>0$) and all positive ones in the lower.) It is clear
that we can deform $D$ into a new diagram $D'$ which represents $L$,
which is such that, when looked upon as a collection of plane curves, it has  
nowhere vanishing curvature, and which satisfies conditions
(a)-(d) of Proposition 1, (i) in ~\cite{V}. Then (ii) of Proposition 1
in ~\cite{V} implies that $D'$ is the diagram of a
holonomic representative of $L$, associated to, say,  $g\colon
C_j\to\tR$. We may assume (see Lemma ~\ref{cod2}) that $g$ is
torsion-generic.  

It will be convenient to look upon $g$ as a function defined on a
collection of intervals rather than on $C_k$. That is,
$g\colon\bigsqcup_{i=1}^j[0,M_i]\to\tR$, 
$g(0)=g(M_i)$ for $1\le i\le j$. 

Let $k>0$ be a real number. Then the collection of holonomic curves
$\left(c(k)\right)(t)=\left(g(kt),kg'(kt),k^2g''(kt)\right)$,
$t\in\bigsqcup_{i=1}^j\left[0,\frac{M}{k}\right]$ also represents $L$. 
As $k\to 0$, $c(k)$ approaches a collection of curves in the
$x_1x_2$-plane with nowhere vanishing curvature and total tangential
degree $2\pi\ba(L)$. It follows that
$$
\lim_{k\to 0}\int_{c(k)}(\kappa + |\tau|)\, ds=2\pi\ba(L).
$$
Thus, $\hct(L)\le2\pi\ba(L)$.\qed     
\section{Proof of Theorem ~\ref{thmhc}}\label{pfthmhc}
The last part of the proof of Theorem ~\ref{thmhct}, implies
that $\hc(L)\le2\pi\ba(L)$. To prove the other inequality 
\begin{equation}\label{low}
A(\ba(L)-\bi(L))+2\pi\bi(L)\le\hc(L) 
\end{equation}
we use the following lemma:
\begin{lma}\label{projmax}
Let $f\colon S^1\to\tR$ have at least $2n$ local
extrema and only non-degenerate critical points. Then for almost every 
$v=(x,y,z)\in \tS^2\subset\tR^3$ such that $y^2-4xz>0$, the function 
$$
h_v(t)=xf(t)+yf'(t)+zf''(t)
$$
has at least $n$ local maxima.
\end{lma}
\begin{pf}
Let $a\in\tR$ be a coordinate on $\tR P^1$, $a\mapsto[1,a]$, where
$[\xi,\eta]$ are homogeneous coordinates on $\tR P^1$. 
Consider the function $g_a(t)=f(t)+af'(t)$. The 
function $f$ has at least $n$ non-degenerate local maxima and $n$
non-degenerate local minima. Between any two local maxima (minima)
there is a  
local minimum (maximum). Thus, if $a\ne 0$ then $g'_a$ changes sign
between the critical points of $f$ and hence $g_a$ has at least
$n$ local maxima (minima). (The latter is of course also true for the
case $a=0$ and  $a=\infty$, $g_{\infty}(t)=f'(t)$.) 

Consider the plane curve $t\mapsto c(t)=\left(f(t),f'(t)\right)$,
$t\in S^1$ (which is regular since the critical points of $f$ are
non-degenerate). If $t\in S^1$ is a critical point of $g_a$ then
the tangent line of $c$ at $c(t)$ is $[-a,1]$. Moreover, $t$
is a degenerate critical point of $g_a$ if and only if the curvature
$\kappa(t)$ of $c$ at $c(t)$ is zero. 

Let $T(t)\in\tR P^1$ denote the tangent line of $c$ at $c(t)$. If
$\kappa^{-1}(0)=U\subset S^1$ then $T(U)$ is exactly the
set of critical values of the map $T\colon S^1\to\tR P^1$. By Sard's
theorem, $T(U)\subset\tR P^1$ is a set of measure zero. Thus, for
almost every $a\in\tR$, $g_a$ has at least $n$ non-degenerate local
maxima (minima). 

Consider the function 
$$
g_{(a,b)}=(f+af')+b(f'+af'')=f+(a+b)\,f'+ab\,f,\quad
(a,b)\in\tR^2.
$$ 
It follows from the above that for almost
every $(a,b)\in\tR^2$, $g_{(a,b)}$ has at least
$n$ local maxima. The map $\tR^2\to\tR P^2$,
$(a,b)\mapsto[1,a+b,ab]$ is a diffeomorphism onto its image when
restricted to $\{(a,b)\colon a>b\}$. The image consists of points
$[1,\eta,\zeta]\in\tR P^2$ such that $\eta^2-4\zeta>0$. This means
that for almost every line $[1,\eta,\zeta]\in\tR P^2$ such that
$\eta^2-4\zeta>0$, the function $f+\eta f'+\zeta f''$ has at least $n$
local maxima (minima).

Taking preimages of the double cover $\tS^2\to\tR P^2$,
$(x,y,z)\mapsto[x,y,z]$ we obtain the statement in the lemma. 
\end{pf}
We now proceed with the proof of ~\eqref{low}: Let
$t\mapsto c(t)=\left(f(t),f'(t),f''(t)\right)$ for some generic
$f\colon C_j\to\tR$ be a holonomic representative of $L$. Let
$C^1,\dots,C^j$ be the components of $C_j$. Let $f_i$ and $c_i$ be the
restrictions of $f$ and $c$, respectively, to $C^i$. Let $n_i$ be the
number of local maxima of $f_i$ and let $n=\sum_{i=1}^jn_i$. Then,
as in the proof of Theorem ~\ref{thmhct}, $n\ge\ba(L)$. 

Since all critical points of a generic function are non-degenerate,
Lemma ~\ref{projmax} implies that for almost every direction 
$v=(x,y,z)$ such that $y^2-4xz>0$, $\mu(c_i,v)\ge n_i$. 
For $v$ in the complementary region ($y^2-4xz\le 0$) such that the
projection onto the 
orthogonal complement of $v$ gives a generic link diagram (a set of
full measure), $\sum_{i=1}^j\mu(c_i,v)\ge\bi(L)$, by definition of
$\bi(L)$. By Lemma ~\ref{avtc}, 
$$
\int_c\kappa\,ds=\frac12\sum_{i=1}^j\int_{\tS^2}\mu(c_i,v)\,dA(v)\ge
A\ba(K)+(2\pi-A)\bi(K). 
$$
\qed

\section{The knots $K(m,j)$}\label{K(m,j)}
\begin{figure}[htbp]
\begin{center}
\includegraphics[angle=0, width=12cm]{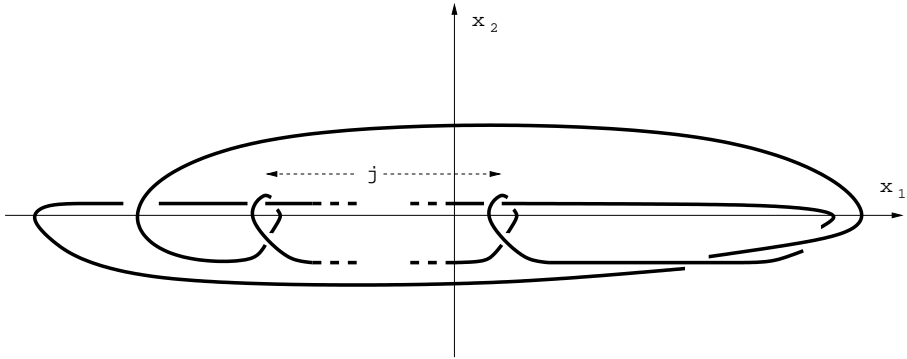}
\end{center}
\caption{Holonomic representatives of the knots $K(2,j)$}\label{kn1} 
\end{figure}

\begin{dfn}
For $j\ge 1$, define the knot $K(2,j)$ to be the class of the
(holonomic) embedding in Figure ~\ref{kn1}.  
\end{dfn}
Note that $K_{2,1}$ is the figure eight knot. 
\begin{dfn}
For $j\ge 1$ and $m\ge 2$, define the knot $K(m,j)$ to be the
connected sum of $K(2,j)$ and $m-2$ right-handed trefoil 
knots (see Figure ~\ref{add}).
\end{dfn}
\section{Proof of Proposition ~\ref{prpex}}\label{pfprpex}
The proof is a combination of lemmas which are stated here and proved
later.   

\begin{lma}\label{babi}
For $m\ge 2$ and $j\ge 1$, the knot $K(m,j)$ satisfies $\bi(K(m,j))=m$ and
$\ba(K(m,j))=m+j$.  
\end{lma}
\noindent
Lemma \ref{babi} is proved in Section \ref{pfbabi}. 

To finish the proof, holonomic representatives of
$K(m,j)$ with total curvature arbitrarily close to $Aj+2\pi m$ must be 
constructed. When this has been done for $m=2$ the 
other cases are easy. We therefore restrict attention to $K(2,j)$.  
Before going into details, the argument will be outlined: 

In Lemma ~\ref{t^n}, we construct a family of non-closed holonomic
curves with the following two properties: First, any curve in the
family can be closed in such way that the result is a closed holonomic 
curve
arbitrarily close to representatives of $K(2,j)$. Second, the
infimum of the total curvatures of the curves in the family is
$Aj+2\pi$. 

The curves in this family are however degenerate at the points where 
they have to be perturbed to give representatives of $K(2,j)$. 
Hence, there is no guarantee that the required perturbation (even if
it is very small) will give a small change in total curvature.  

In Lemma ~\ref{t^3}, we overcome this problem: By deforming the
above curves in a specific way, a new family of curves is
produced. The new family have the properties of the old family and, in
addition to that, its members are curves which are non-degenerate
where they must be perturbed to give representatives of $K(2,j)$. This
non-degeneracy warrants that the total curvature is 
continuous under small perturbations.

Finally, we pick a curve in the new family with total curvature close
enough to $Aj+2\pi$ and make a small enough perturbation. In Lemma
~\ref{join} we prove that it is possible to close the perturbed curve
in such a way that the result is a holonomic representative of
$K(2,j)$ with total curvature as close as required to $Aj+4\pi$.

We now go into details: 
Let $M\ge 0$ and let $a_1,\dots,a_j$ be points in
$\left[0,M\right]$.  
For $i=1,\dots,j$, let $J_i=\left[a_i-2,a_i+2\right]$ and
$I_i=\left[a_i-1,a_i+1\right]$. Assume that the 
$J_i\cap J_l=\emptyset$ for $i\ne l$ and that
$J_i\subset\left(0,M\right)$ for every $i$. 

\begin{figure}[htbp]
\begin{center}
\includegraphics[angle=0, width=12cm]{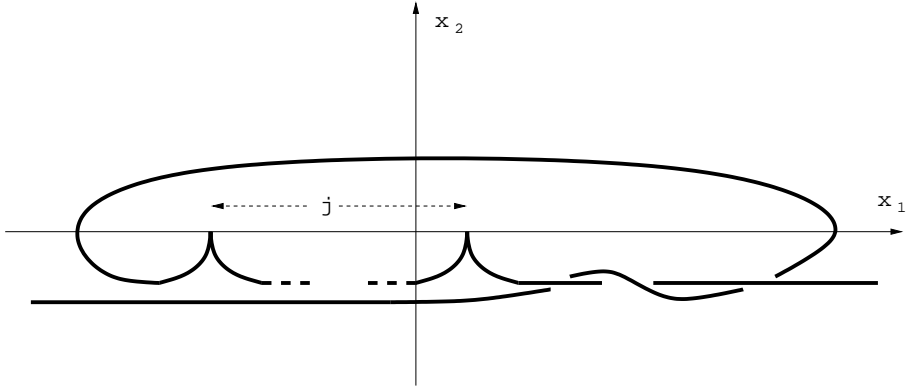}
\end{center}
\caption{The holonomic curve associated to $f_{j,n}$}\label{cur} 
\end{figure}

For {\em odd} $n\ge 3$, let $f_{j,n}\colon \left[0,M\right]\to\tR$ be
a function which satisfies 
$$
f_{j,n}(t)=-(t-a_i)^n \quad\text{ for } t\in I_i,\quad 1\le i\le j, 
$$
which has constant derivative in neighborhoods of the endpoints of 
$\left[0,M\right]$, and which has associated 
holonomic curve is as in Figure ~\ref{cur}. (Note that the associated
holonomic curve is non-closed).  

For $k>0$, let $c(j,n,k)$ be the holonomic curve associated to the
function $t\mapsto f_{j,n}(kt)$, $t\in \left[0,\frac{M}{k}\right]$. 
Note that for $n>3$, the curvature function of $c(j,k,n)$ is not
defined at the points $\frac{a_i}{k}$, where $c(j,n,k)$ has contact of
order $n-2$ with a line in the $x_3$-direction. 
\begin{lma}\label{t^n}
$$
\lim_{k\to 0}\int_{c(j,n,k)}\kappa\,ds= A(n)j+2\pi, 
$$
where $2A(n)$ is the area on the unit sphere of the region 
$\{(x,y,z)\in\tS^2\colon y^2-4\left(\frac{n-2}{n-1}\right)xz>0\}$. 
\end{lma}
\noindent
Lemma ~\ref{t^n} is proved in Section ~\ref{pft^n}.

Let $p\colon[0,M]\to\tR$ be a function which satisfies 
$$
p(t)=\begin{cases}
      -(t-a_i)^3 & \text{for } t\in I_i,\quad 1\le i\le j,\\ 
         0       & \text{for } t\in \left[0,M\right]-\bigcup_{i=1}^jJ_i.
      \end{cases}
$$ 
For $\delta>0$, let $c(j,n,k,\delta)$ denote the
holonomic curve associated to the function  
$t\mapsto f_{j,n}(kt)+\delta p(kt)$, $t\in\left[0,\frac{M}{k}\right]$. 
Clearly, the curve $c(j,n,k,\delta)$ is non-degenerate on 
$\left[\frac{a_i-1}{k},\frac{a_i+1}{k}\right]$,
$1\le i\le j$, if $\delta>0$ is small enough.  

\begin{lma}\label{t^3}
Let $0<a<1$. If $n\ge 2\left(1-\sqrt{1-2a^2}\right)^{-1}+2$ and 
$0\le k\le 1$  
then  
$$
\lim_{\delta\to 0}\int_{c(j,n,k,\delta)}\kappa\,ds\le
\int_{c(j,n,k)}\kappa\,ds + 6B(a)j,
$$
where $2B(a)$ 
is the area on the unit sphere of the region
$
\{(x,y,z)\in \tS^2\colon |x|\le a \text{\rm{ or }} |y|\le a 
\text{\rm{ or }} |z|\le a\}
$. 
\end{lma} 

\noindent
Lemma ~\ref{t^3} is proved in Section ~\ref{pft^3}.

Let $q\colon\left[0,M\right]\to\tR$ be a function which satisfies
$$
q(t)=\begin{cases}
     (t-a_i) & \text{for } t\in I_i,\quad 1\le i\le j,\\ 
        0    & \text{for } t\in \left[0,M\right]-\bigcup_{i=1}^jJ_i.
     \end{cases} 
$$
For $\beta>0$, let 
$c(j,n,k,\delta,\beta)$ denote the holonomic curve associated to the
function $t\mapsto f(kt)+\delta p(kt)+\beta q(kt)$,
$t\in\left[0,\frac{M}{k}\right]$. 
Fix a small $\delta>0$ such that $c(j,n,k,\delta)$
is non-degenerate along the support of $t\mapsto p(kt)$. This
non-degeneracy implies that the total curvature is continuous under
small perturbation and therefore: 
\begin{equation}\label{limbeta}
\lim_{\beta\to 0}\int_{c(j,n,k,\delta,\beta)}\kappa\,ds=
\int_{c(j,n,k,\delta)}\kappa\,ds.
\end{equation} 

\begin{figure}[htbp]
\begin{center}
\includegraphics[angle=0, width=8cm]{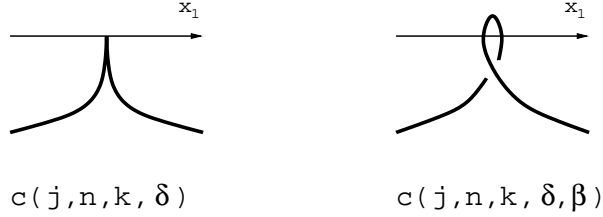}
\end{center}
\caption{The effect of adding $\beta q$ to $f_{j,n}+\delta p$}\label{beta} 
\end{figure}

\begin{figure}[htbp]
\begin{center}
\includegraphics[angle=0, width=12cm]{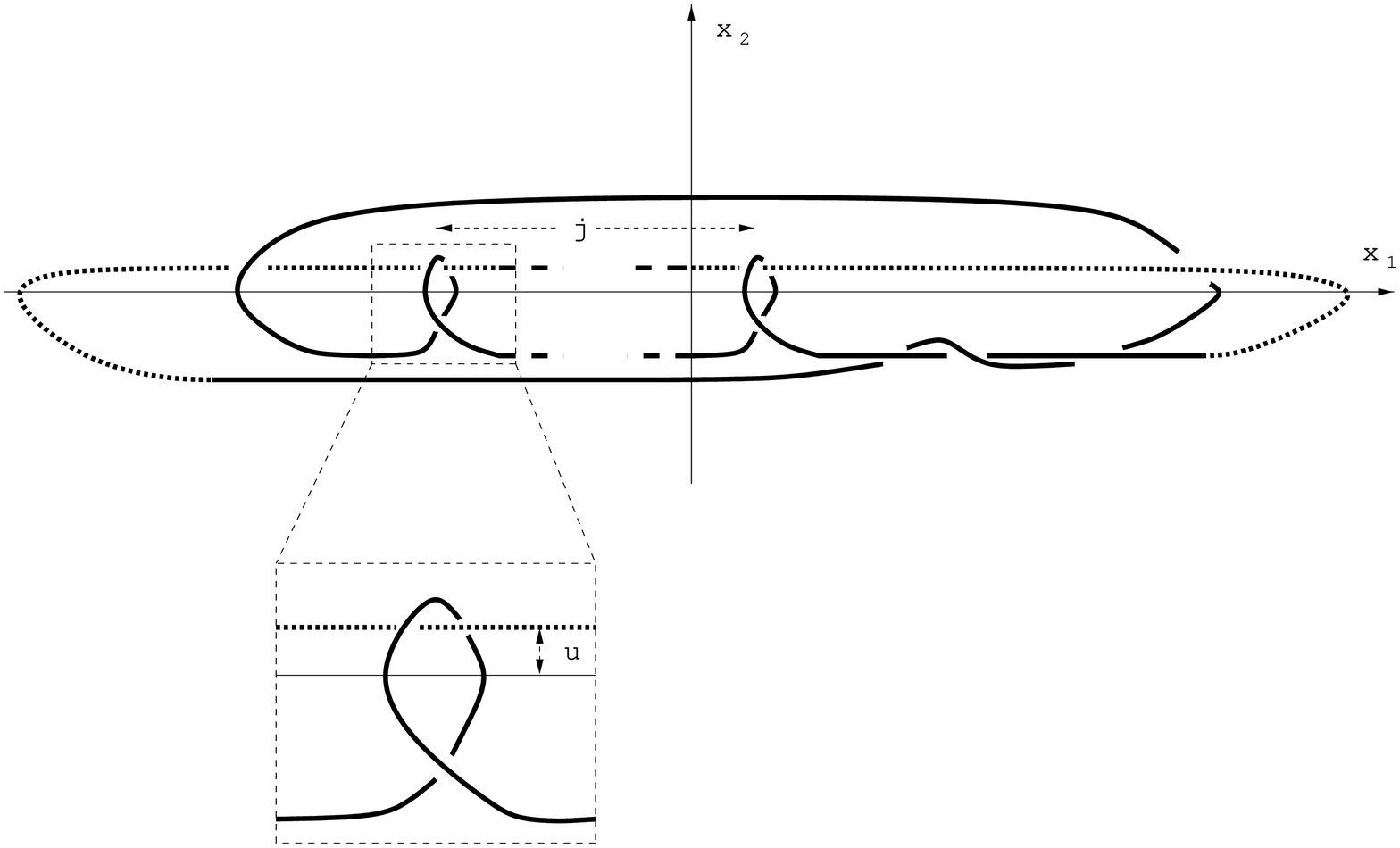}
\end{center}
\caption{Closing the curve}\label{clo} 
\end{figure}

As mentioned in the outline above, we need to close the holonomic
curve\linebreak 
$c(j,n,k,\delta,\beta)$ so that the result is a holonomic
representative of the knot $K(2,j)$. It is easy to see that this is
possible (see Figures ~\ref{beta} and 
~\ref{clo}). The next lemma allows us to control the total curvature
when closing the curve. 

\begin{lma}\label{join}
Let $c,u,d>0$ be given constants and let $s,\alpha>0$. 
Then there exist $S>0$ and functions $f,g\colon\left[-s,S+s\right]\to\tR$ such
that: 
\begin{itemize}
\item[{\rm (i)}] For $t\in[-s,0]$, $f(t)=ut+c$ and  $g(t)=-dt-c$. 
\item[{\rm (ii)}] For $t\in\left[S,S+s\right]$, $f(t)=-d(t-S)+c$ and
$g(t)=u(t-S)-c$. 
\item[{\rm (iii)}] For all $t$, $|f'(t)|\le\max(u,d)$ and
$|g'(t)|\le\max(u,d)$. 
\item[{\rm (iv)}] The total curvatures of the holonomic curves
associated to $f$ and $g$ respectively are smaller than $\pi+\alpha$.
\end{itemize}
\end{lma}
\noindent
Lemma \ref{join} is proved in Section \ref{pfjoin}.

We are now in position to prove Proposition ~\ref{prpex}: 
The statements about bridge- and braid-index is just Lemma
~\ref{babi}. 

Let
$\epsilon>0$ be given. We must find a 
holonomic representative of $K(m,j)$ with total curvature less
than $Aj+2\pi m+\epsilon$. We start with the case $K(2,j)$: 

Choose $a>0$ such that $B(a)<\frac{\epsilon}{48j}$. Choose 
$n\ge 2\left(1-\sqrt{1-2a^2}\right)^{-1}+2$ such that
$A(n)-A\le\frac{\epsilon}{8j}$. By Lemma 
~\ref{t^n} it is possible to find $k>0$ such that
$$
\int_{c(j,n,k)}\kappa\,ds\le Aj+2\pi+\frac{\epsilon}{4}.
$$ 
Lemma ~\ref{t^3} then implies that for $\delta>0$ small enough
$$ 
\int_{c(k,n,\delta)}\kappa\,ds\le Aj+2\pi+\frac{\epsilon}{2},
$$
and $c(j,n,k,\delta)$ is non-degenerate along the support of 
$t\mapsto q(kt)$, $t\in\left[0,\frac{M}{k}\right]$. By equation
~\eqref{limbeta}, for $\beta>0$ small enough
$$
\int_{c(j,n,k,\delta,\beta)}\kappa\,ds\le Aj+2\pi+\frac{3\epsilon}{4}.
$$
Finally, we close $c(k,n,\delta,\beta)$ as in Figure
~\ref{clo} (we add the dashed part). 
If $u$ in Figure ~\ref{clo} is
small enough then the resulting holonomic curve is a representative of
$K(2,j)$. Using Lemma ~\ref{join}, we can assure that it has 
total curvature less than $Aj+4\pi+\epsilon$. This proves the
Proposition for $K(2,j)$. 

\begin{figure}[htbp]
\begin{center}
\includegraphics[angle=0, width=12cm]{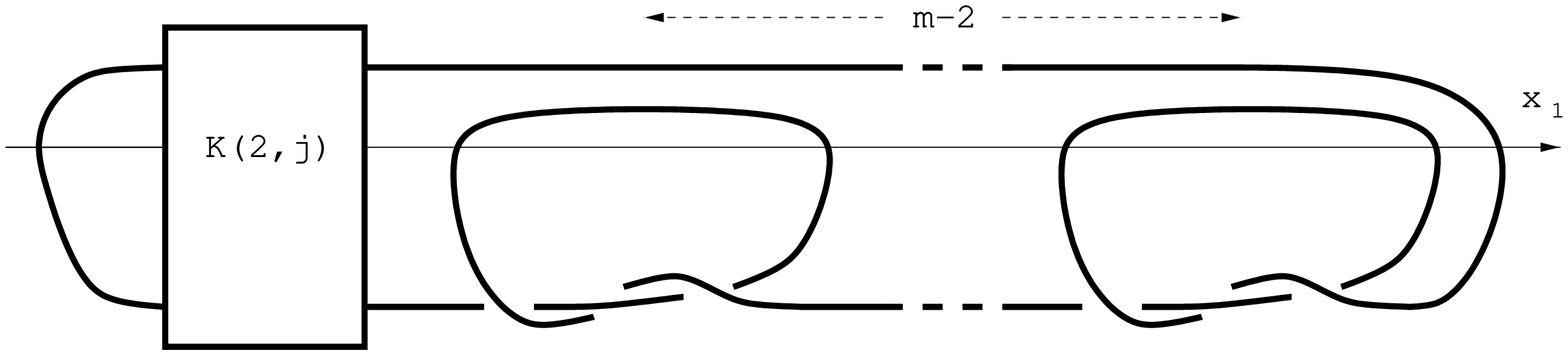}
\end{center}
\caption{Adding trefoils to $K(2,j)$}\label{add} 
\end{figure}

The general case follows by noting that it is possible to add $m-2$
holonomic trefoils to a holonomic representative of $K(2,j)$, with
total curvature smaller than $Aj+4\pi+\frac{\epsilon}{2}$, in such a 
way that the increase in total curvature is smaller
than $2\pi(m-2)+\frac{\epsilon}{2}$ (see Figure ~\ref{add}). \qed 
\subsection{Proof of Lemma \ref{babi}}\label{pfbabi}
Consider the representative of $K(2,j)$ presented in Figure
~\ref{kn1}. If $j$ third Reidemeister moves, moving all crossings
except the three rightmost into the upper half plane are followed by
$j$ first Reidemeister moves removing the $j$ uppermost crossings
then the resulting (non-holonomic) embedding has exactly two maxima in
the $x_1$-direction. It follows that $\bi(K(2,j))\le 2$ but $K(2,j)$
is knotted. Thus, $\bi(K(2,j))=2$.  

The bridge-index of the trefoil knot is $2$ and hence $\bi(K(m,j))=m$, 
by the additivity properties of bridge-indices.  

Looking at the projection of the representative of $K(2,j)$ in
Figure ~\ref{kn1} onto the $x_2x_3$-plane, we see that
$\ba(K(m,j))\le m+j$. 

A straightforward induction shows that the HOMFLY-polynomial
$P_{K(2,j)}(z,v)$ of $K(2,j)$, written as a Laurent 
polynomial in $v$ with coefficients in $\tZ[z,z^{-1}]$, 
is of the form
$$
P_{K(2,j)}(z,v)=v^{-2j}+\dots+v^2,
$$
where $\dots$ indicate a Laurent polynomial in $v$ with terms of
degrees strictly between the degrees of the terms written out. The
HOMFLY-polynomial of the right-handed trefoil is $(2+z^2)v^2-v^4$. 
Since the HOMFLY-polynomial is multiplicative under connected sum,
$$
P_{K(m,j)}=(2+z^2)^{m-2}\, v^{2(m-j-2)}+\dots
+(-1)^{m-2}\, v^{2(2m-3))}. 
$$
Now, for any knot $K$, $\ba(K)\ge\frac12\text{$v$-span}(P_K(z,v))+1$ 
(see Morton ~\cite{Mo}). Hence, $\ba(K(m,j))\ge m+j$.\qed
\subsection{Proof of Lemma ~\ref{t^n}}\label{pft^n}
Fix $n$. Let $U(k)$ denote the union of intervals 
$$
U(k)=\left[0,\frac{M}{k}\right]-
\bigcup_{i=1}^j
\left[\frac{a_i-1}{k},\frac{a_i+1}{k}\right].
$$ 
As $k\to 0$ the restriction of the curve $c(j,n,k)$ to $U(k)$,
approaches parts of a plane curve with all its curvature concentrated
at intersections with the $x_1$-axis. In the limit, this part
of the curve contributes $2\pi$ to the total curvature (see Figure
~\ref{cur}).  

On $\left[\frac{a_i-1}{k},\frac{a_i+1}{k}\right]$,
$c(j,n,k)$ is a translate of the holonomic curve associated to
$t\mapsto -k^nt^n$ on $\left[-\frac{1}{k},\frac{1}{k}\right]$. It is
straightforward to see that the factor $-k^n$ does not affect the
total curvature of this curve. Hence, in the limit as $k\to 0$,
the contribution to the total curvature from each of these parts of
$c(j,n,k)$ is
\begin{equation}\label{b(n)}
\int_{b(n)}\kappa\,ds,
\end{equation}  
where $b(n)$ is the curve $(t^n,nt^{n-1},n(n-1)t^{n-2})$,
$t\in\tR$. 

We use Lemma ~\ref{avtc} to evaluate the integral ~\eqref{b(n)}.
For a unit vector $v=(x,y,z)$, let  
$$
h_v(t)=xt^n+nyt^{n-1}+n(n-1)zt^{n-2}, 
$$
It is straightforward to check that $h_v$ has one local maximum for
almost every $v$ such that $y^2-4\left(\frac{n-2}{n-1}\right) xz>0$ and
no local maximum for $v$ in the complementary region on $\tS^2$. \qed 
\subsection{Proof of Lemma ~\ref{t^3}}\label{pft^3}
Assume $n\ge 5$. Let
$$
J=\bigcup_{i=1}^j
\left[\frac{a_i-2}{k},\frac{a_i+2}{k}\right],\quad
I=\bigcup_{i=1}^j
\left[\frac{a_i-1}{k},\frac{a_i+1}{k}\right],\quad\text{and }
U=\left[0,\frac{M}{k}\right]-I
$$ 
The function $t\to p(kt)$, $t\in\left[0,\frac{M}{k}\right]$ has
support in $J$ and the curve $c(j,n,k)$ can be taken non-degenerate in
$J-I$. Therefore,  
\begin{equation}\label{Upart}
\lim_{\delta\to 0}\int_{c(j,n,k,\delta)|U}\kappa\,ds=
\int_{c(j,n,k)|U}\kappa\,ds,
\end{equation}
where $|U$ denotes restriction to $U$. 

We are left with the image of $I$. On the intervals in $I$, 
$t\mapsto f_n(kt)+\delta p(kt)$ is a translate of   
the function $t\mapsto -k^nt^n-\delta t^3$,
$t\in\left[-\frac{1}{k},\frac{1}{k}\right]$. Using Lemma ~\ref{avtc} as in
the proof of Lemma ~\ref{t^n} we need to calculate the number of
maxima of the function  
$$
h_v(t)=k^nnt^{n-2}\left(xt^2+(n-1)yt+(n-1)(n-2)z\right)+
\delta t\left(xt^2+2yt+2z\right), 
$$ 
for $v=(x,y,z)\in\tS^2$. Taking derivatives and assuming 
$x\ne 0$ we must calculate the number of zeros of the polynomial
\begin{align*}
h'_v(t)&=k^nnt^{n-3}
\left(t^2+(n-1)\frac{y}{x}t+(n-1)(n-2)\frac{z}{x}\right)+
3\delta\left(t^2+2\frac{y}{x}t+2\frac{z}{x}\right)\\
&=k^nnt^{n-3}g(t)+3\delta r(t),
\end{align*}
where the last equality is used to define the quadratic polynomials
$g$ and $r$.

Below we make several observations which together give an estimate on 
the total curvature of the image of $I$:

{\em Observation 1}: For any $v$, $h'_v$ cannot have more than $6$ real
zeros. This follows from Descartes' Lemma (see for example Benedetti
and Risler ~\cite{BR}, Proposition 1.1.10).    

{\em Observation 2}:
If $y^2-2xz<0$ then neither $g$ nor $r$ have any real zeros. Hence,
$h'_v$ does not have real zeros in this case. 

{\em Observation 3}:
If $y^2-4\left(\frac{n-2}{n-1}\right)xz>0$ and $xz<0$ then both $g$
and $r$ have two real zeros on opposite sides of $t=0$. Moreover, the
zeros of $r$ 
are closer to $t=0$ than those of $g$ and $t\mapsto t^{n-3}g(t)$ has a 
maximum at $t=0$. It follows that $h'_v$ has exactly two real
zeros in this case. 

{\em Observation 4}:
If $y^2-4\left(\frac{n-2}{n-1}\right)xz<0$ and $y^2-2xz<0$ then
$xz>0$, $g$ has no real zeros, and $r$ has two real zeros on the same
side of $t=0$. Denote the zeros of $r$ by $\theta_1\le \theta_2$.

{\em Observation 5}:
If $y^2-4\left(\frac{n-2}{n-1}\right)xz<0$ and $xz>0$ then 
$t\mapsto t^{n-3}g(t)$ has one local 
maximum at $t=\phi$, $r$ has two real zeros $\theta_1\le \theta_2$,
and $\phi$, $\theta_1$ and $\theta_2$ all lie on the same side of
$t=0$. 

If $v$ satisfies $|x|>a$, $|y|>a$, and $|z|>a$, and if $n$ is large
enough then $|\phi|>\max(|\theta_1|,|\theta_2|)$: The derivative of 
$t\mapsto t^{n-3}g(t)$ is 
$$
\frac{d}{dt}\left(t^{n-3}g(t)\right)=
(n-1)t^{n-4}\left(t^2+(n-2)\frac{y}{x}t+(n-2)(n-3)\frac{z}{x}\right),
$$ 
which has zeros at 
$$
t=0\quad\text{and}\quad t=\frac{n-2}{2}
\left(-\frac{y}{x}
\pm\frac{1}{|x|}\sqrt{y^2-4\left(\frac{n-3}{n-2}\right)xz}\right).
$$
It is straightforward to check that the non-zero zeros have distance
at least\linebreak 
$\frac{n-2}{2|x|}\left(1-\sqrt{1-2a^2}\right)$ from  
$t=0$. The zeros of $r$ are
$t=-\frac{y}{x}\pm\frac{1}{|x|}\sqrt{y^2-2xz}$ 
which have distance at most $\frac{2}{|x|}$ from $t=0$. Thus, if 
$$
n>2\left(1-\sqrt{1-2a^2}\right)^{-1}+2
$$  
then $|\phi|>\max(|\theta_1|,|\theta_2|)$.

{\em Observation 6}:
If $\delta>0$ is small enough and $v$
satisfies $|x|>a$, $|y|>a$, and $|z|>a$ then, with $v$ as in
Observation 4 or Observation 5, $h'_v$ is monotonic on
$\left[\theta_1,\theta_2\right]$ (which implies that $h'_v$ has no zeros
if $v$ is as in Observation 4 and no more than two zeros if $v$ is as
in Observation 5): 

The derivative $h''_v$ of $h'_v$ is
\begin{align*}
h''_v(t)&=
k^nn(n-1)t^{n-4}\left(t^2+(n-2)\frac{y}{x}t+(n-2)(n-3)\frac{z}{x}\right) 
+6\delta\left(t+\frac{y}{x}\right)\\
&=l(t)+3\delta r'(t),
\end{align*}
where the last equality serves as a definition of $l$.
We must check that this expression does not change sign on
$\left[\theta_1,\theta_2\right]$. On this interval
$|r'(t)|<\frac{2}{a}\sqrt{1-2a^2}$ and 
$|l(t)|\ge |l\left(1-\sqrt{1-2a^2}\right)|>0$. 
It follows that $h'_v$ is monotonic
on $\left[\theta_1,\theta_2\right]$ for $\delta$ small enough.  

We collect the observations above to a proof:
Assume $v$ is such that $|x|>a$, $|y|>a$, and $|z|>a$. Then 
Observations  2, 4, and 6 imply that if
$y^2-4\left(\frac{n-1}{n-2}\right)xz<0$ then  
$h_v$ does not have any maxima, and Observations 3, 5, and 6 imply that
if $y^2-4\left(\frac{n-1}{n-2}\right)xz>0$ then $h_v$ has at most
one maximum. Assume that $v$ is such that $|x|\le a$, or $|y|\le a$, or
$|z|\le a$. Then Observation 1 implies that $h_v$ has at most $6$
maxima. 
Thus, by Lemma ~\ref{avtc},
\begin{equation}\label{Ipart}
\lim_{\delta\to 0}\int_{c(n,k,\delta)|I}\kappa\,ds\le
\int_{c(n,k)|I}\kappa\,ds+6B(a)j.
\end{equation}
The lemma follows from ~\eqref{Upart} and ~\eqref{Ipart}.\qed     
\subsection{Proof of Lemma ~\ref{join}}\label{pfjoin}
We construct the function $f$ (the function $g$ can be constructed in
a similar way):

Let $(x_1,x_2,x_3)$ be coordinates on $\tR^3$. To fulfill condition (iv), 
we must find $f$ such that the length of the curve $e_1(t)$, 
$0\le t\le S$ traced out by the unit tangent vector of
$t\mapsto c(t)=\left(f(t),f'(t),f''(t)\right)$ on $\tS^2$ is arbitrary 
close to $\pi$. 

Note that $c'(t)=\left(f'(t),f''(t),f'''(t)\right)$. Condition (i)
implies that $e_1(0)=(1,0,0)$ and condition (ii) that $e_1(S)=(0,0,-1)$.  

We first make a non-smooth model: Let $w=\min(u,d)$ and let $K\ge 1$
(ultimately we shall take $K$ very large). Let
$S=\sqrt{\frac{2K}{uw}}(u+d)$ and let  
$$
h(t)=
\begin{cases}
-\frac{w}{K} & \text{for } 0\le t\le \sqrt{\frac{2k}{uw}}u\\
\frac{uw}{dK} & \text{for } \sqrt{\frac{2k}{uw}}u< t\le
\sqrt{\frac{2k}{uw}}(u+d). 
\end{cases}
$$ 
Take $f'''(t)=h(t)$ on $\left[0,S\right]$ and $f'''(t)=0$ for other
$t$. With $f''(0)=0$ and $f'(0)=u$ this defines $f'$ on
$\left[-s,S+s\right]$. It is easy to check that 
\begin{equation}\label{band}
\frac{|f'''|}{\sqrt{(f')^2+(f'')^2}}\le\frac{8w}{d\sqrt{K}},
\end{equation}    
for all $t\in[-s,S+s]$.
Inequality ~\eqref{band} implies that the
curve $e_1$ lies entirely inside a  band 
around the equator of $\tS^2$ in the $x_1x_2$-plane of height (in
the $x_3$-direction) less than $\frac{8w}{d\sqrt{K}}$. It follows that
the length of the curve $e_1$ approaches $\pi$ as $K\to\infty$. 

It is clearly possible to smooth the above model keeping the total
curvature close to $\pi$. 

This construction gives the derivative of the desired function. Denote
this derivative $k(t)$. If $\int_0^Sk(t)\,dt=r>0$ then we take
$S'=S+\frac{r}{d}$ and let $f(t)=c+\int_0^Sk(t)\,dt$ on $[0,S]$ and 
$f(t)=-d(t-S)+r$ on $\left[S,S'\right]$ and we get a function
satisfying (a)-(d) on $\left[-s,S'+s\right]$. If
$\int_0^Sk(t)\,dt=r<0$ we can proceed in a similar way, changing the
function close to $t=0$.\qed
\section*{Acknowledgments}
Parts of this paper is the master thesis of the second author. The
first author would like to thank Hitoshi Murakami for valuable 
discussions and comments. During the preparation of this paper the
first author was supported by the Royal Swedish Academy of Sciences.

\end{document}